\documentclass[a4paper,11pt,reqno]{amsart}
\usepackage{amssymb,amsmath,hyperref}
\usepackage{enumerate,enumitem}
\usepackage{multicol}
\usepackage{cleveref}
\usepackage{tikz,tkz-euclide, circuitikz}
\usepackage{float}
\usetikzlibrary{arrows,arrows.meta,patterns,matrix,
	shapes,fit,calc,shadows,plotmarks}
\author{Anton Bugleev}
\address{Department of Mathematics and Computer Science, Saint Petersburg State University, Saint Petersburg,  Russia, 199178}
\email{handelluss@gmail.com}

\renewcommand\theta{\vartheta}

\newtheorem{theorem}{Theorem}

\theoremstyle{definition}

\newtheorem{remark}[theorem]{Remark} 
\newtheorem{example}[theorem]{Example}
\newtheorem{conjecture}[theorem]{Conjecture}

\numberwithin{theorem}{section} 
\numberwithin{equation}{section}

\crefname{part}{Part}{Parts}
\Crefname{part}{Part}{Parts}

\makeatletter
\@namedef{subjclassname@2020}{\textup{2020} Mathematics Subject Classification}
\makeatother
\begin{document}

\title{Two-color partitions and overpartitions: a combinatorial proof}
\date{9 September 2025}


\tikzset{
    up triangle/.pic={
        \draw (0,0) -- (1,0) -- (1,-1) -- cycle;
        \node at (0.7,-0.3) {};
    },
    down triangle/.pic={
        \draw (0,0) -- (0, -1) -- (1, -1) -- cycle;
        \node at (0.3,-0.7) {};
    },
    square/.pic={
        \draw (0,0) -- (1,0) -- (1,-1) -- (0, -1) -- cycle;
        \node at (0.5,-0.5) {};
    }
}

\begin{abstract}
George Andrews and Mohamed El Bachraoui recently explored identities for two-color partitions \cite{andrews2024twocolorpartitionsoddsmallest}. In particular, they studied the connection between two-colored partitions and overpartitions. Their proofs were analytical, but they also conjectured combinatorial proofs of their results. In this paper we use two-modular diagrams to give a combinatorial proof of their main result.
\end{abstract}

\maketitle
\tableofcontents 

\section{Introduction}

Throughout the paper we let \(n\) denote a nonnegative integer. Recall that {\it a partition} of \(n\) is a tuple \((\lambda_1, \dots, \lambda_k)\) such that \(\lambda_1 + \dots + \lambda_k = n\) and \(\lambda_1 \geq \lambda_2 \geq \dots \geq \lambda_k \geq 1\). Let us call a partition {\it two-colored} if each summand is colored in either green or blue. For example, the tuple \((9_g, 5_b, 4_g, 2_g, 2_b, 2_b)\) is a partition of 24.

{\it An overpartition} of \(n\) is a partition in which first occurence of each number may be overlined. For example, \((5, \overline{3}, 3, 2, 1, 1)\) is an overpartition of 15.

Let \(n\) be a nonnegative integer and let \(\mathcal{E}(n)\) be the set of all two-color partitions of \(n\) in which parts are distinct and the even parts may occur only in the blue color. Let \(\overline{\mathcal{P}_o}(n)\) be the set of all overpartitions of \(n\) into odd parts. Denote the cardinality of \(\mathcal{E}(n)\) and \(\overline{\mathcal{P}_o}(n)\) by \(E(n)\) and \(\overline{p_o}(n)\) respectively.

Let \(E_0(n)\) (respectively \(E_1(n)\)) denote the number of partitions counted by \(E(n)\) in which the number of even parts is even (respectively odd). Let \(E_2(n)\) (respectively \(E_3(n)\)) denote the number of partitions counted by \(E(n)\) in which the number of all parts is even (respectively odd).

\begin{remark}
    Throughout this paper a two-color partition has distinct parts if and only if any two parts of the partition have different colors or if they are numerically distinct.
\end{remark}

Recently George E. Andrews and Mohamed El Bachraoui \cite{andrews2024twocolorpartitionsoddsmallest} established the following connection between two-colored partitions and overpartitions. 
\newpage
\begin{theorem} \label{thm}
    For any nonnegative integer \(n\) the following holds:
    \begin{enumerate}[label=\alph*)]
    \item \label[part]{thm:a}\(E(n) = \overline{p_o}(n)\),
    \begin{multicols}{2}
        \item \label[part]{thm:b} \(E_0(n) = 
        \begin{cases}
            \frac{\overline{p_o}(n)}{2} + 1 & \text{if \(n\) is a square}, \\
            \frac{\overline{p_o}(n)}{2} & \text{otherwise},
        \end{cases}\)
        \item \label[part]{thm:c}\(E_1(n) = 
        \begin{cases}
            \frac{\overline{p_o}(n)}{2} - 1 & \text{if \(n\) is a square}, \\
            \frac{\overline{p_o}(n)}{2} & \text{otherwise},
        \end{cases}\)
        \item \label[part]{thm:d}\(E_2(n) = 
        \begin{cases}
            \frac{\overline{p_o}(n)}{2} + (-1)^n & \text{if \(n\) is a square}, \\
            \frac{\overline{p_o}(n)}{2} & \text{otherwise},
        \end{cases}\)
        \item \label[part]{thm:e}\(E_3(n) = 
        \begin{cases}
            \frac{\overline{p_o}(n)}{2} - (-1)^n &  \text{if \(n\) is a square},\\
            \frac{\overline{p_o}(n)}{2} & \text{otherwise}.
        \end{cases}\)
    \end{multicols}
    \end{enumerate}
\end{theorem}

Their proof is analytical; however, they conjectured a combinatorial proof of \Cref{thm}. In particular, they conjectured
\begin{conjecture}
    (Section 9, \cite{andrews2024twocolorpartitionsoddsmallest}). It would be interesting to give combinatorial proofs for the identities in \Cref{thm}.
\end{conjecture}
In \Cref{sec:2}, we establish a bijection to provide a combinatorial proof of \Cref{thm:a} of \Cref{thm}. In \Cref{sec:3}, we introduce two-modular diagrams. In \Cref{sec:4}, we provide a combinatorial proof of \Crefrange{thm:b}{thm:e} of \Cref{thm}.

\section{Combinatorial proof of \cref{thm:a}} \label{sec:2}

\begin{proof}
Fix any nonnegative integer \(n\) and any overpartition \(\lambda \in \overline{\mathcal{P}_o}(n)\). Recall that an overpartition is a partition in which the first occurrence of any number may be overlined. Let the overlined parts of \(\lambda\) be in the color green and the non-overlined parts of \(\lambda\) be in the color blue. Then one can say that \(\lambda = (\lambda_g, \lambda_b)\) is a two-color partition in which blue parts are odd and green parts are odd and numerically distinct (here \(\lambda_g\) and \(\lambda_b\) denote green and blue parts of \(\lambda\) respectively). 

Recall that Euler's theorem \cite{andrews} states that the number of partitions of a nonnegative integer \(n\) into odd parts coincides with the number of partitions of \(n\) into distinct parts. This theorem has a well-known combinatorial proof.

Then one can say that by Euler's theorem one can map odd blue parts \(\lambda_b\) bijectively to numerically distinct blue parts \(\lambda_b'\). Hence \(\lambda' = (\lambda_g, \lambda_b')\) is a two-color partition into distinct blue parts and odd distinct green parts. In other words \(\lambda' = (\lambda_g, \lambda_b') \in \mathcal{E}(n)\). Clearly the mapping \(\lambda \mapsto \lambda'\) is a bijection.
\end{proof}

We point out that D. Chen and Z. Zou \cite{articlecombinatoricalproofs} also gave a combinatorial proof of \cref{thm:a}.

\section{Two-modular diagrams} \label{sec:3}

Fix any \(\lambda \in \mathcal{E}(n)\). Let us say that \(\lambda = (\lambda_e, \lambda_g, \lambda_b)\), where \(\lambda_e, \lambda_g \text{ and } \lambda_b\) denote blue even parts, green odd parts and blue odd parts of \(\lambda\) respectively.

Let us introduce the {\it two-modular diagram} of \(\lambda\). In these diagrams \(\lambda_g\) parts are represented in rows above the main diagonal and \(\lambda_b\) parts are represented in columns below the main diagonal (\(\lambda_e\) parts are not represented in the diagram). Each odd part is represented in the diagram the following way: if a part equals \(2k+1\), then its representation would have \(k\) squares and 1 triangle. In other words, a square in the diagram represents the number \textquotedblleft {two}\textquotedblright \ and a triangle represents the number \textquotedblleft one\textquotedblright.

\newpage
\begin{example}
\(\lambda_e = \{12, 6\}; \lambda_g = \{13, 9, 5, 3\}; \lambda_b = \{5, 1\}\). 
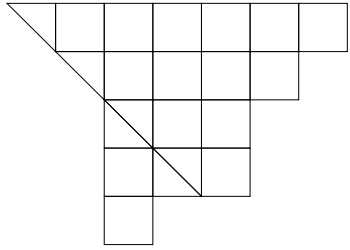
\begin{figure}[H]
\resizebox{130pt}{!}{
\begin{tikzpicture}
\foreach \x in {0, 1, 2, 3} {
    \pic at (\x, 3 - \x) {up triangle};
}
\foreach \x in {1, 2, 3, 4, 5, 6} {
    \pic at (\x,3) {square};
}
\foreach \x in {2, 3, 4, 5} {
    \pic at (\x,2) {square};
}
\foreach \x in {3, 4} {
    \pic at (\x,1) {square};
}
\pic at (4,0) {square};
\pic at (3,0) {down triangle};
\pic at (2,1) {down triangle};
\pic at (2,0) {square};
\pic at (2,-1) {square};
\end{tikzpicture}}    
\end{figure}
\end{example}

\begin{example}
\(\lambda_e = \{12, 4, 2\}; \lambda_g = \{7, 3\}; \lambda_b = \{11, 7, 5, 3\}\).
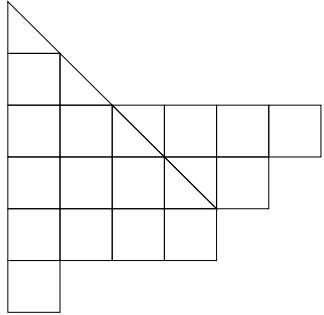
\begin{figure}[H]
\resizebox{120pt}{!}{
\begin{tikzpicture}
\foreach \x in {0, 1, 2, 3} {
    \pic at (\x, 3 - \x) {down triangle};
}
\pic at (2, 1) {up triangle};
\pic at (3, 0) {up triangle};
\foreach \y in {-2, -1, 0, 1, 2} {
    \pic at (0, \y) {square};
}
\foreach \y in {-1, 0, 1} {
    \pic at (1, \y) {square};
}
\foreach \y in {-1, 0} {
    \pic at (2, \y) {square};
}
\pic at (3, -1) {square};
\pic at (4, 0) {square};
\pic at (3, 1) {square};
\pic at (4, 1) {square};
\pic at (5, 1) {square};
\end{tikzpicture}}    
\end{figure}
\end{example}

\begin{example}
\(\lambda_e = \{4, 2\}; \lambda_g = \{11, 7, 3, 1\}; \lambda_b = \{9, 7, 3, 1\}\).
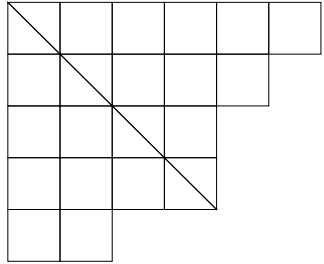
\begin{figure}[H]
\resizebox{120pt}{!}{
\begin{tikzpicture}
\foreach \x in {0, 1, 2, 3} {
\pic at (\x, 3 - \x) {down triangle};
\pic at (\x, 3 - \x) {up triangle};
}
\foreach \x in {1, 2, 3, 4, 5} {
\pic at (\x, 3) {square};
}
\foreach \x in {2, 3, 4} {
\pic at (\x, 2) {square};
}
\pic at (3, 1) {square};
\foreach \y in {-1, 0, 1, 2} {
\pic at (0, \y) {square};
}
\foreach \y in {-1, 0, 1} {
\pic at (1, \y) {square};
}
\pic at (2, 0) {square};
\end{tikzpicture}}    
\end{figure}
\end{example}

\section{Combinatorial proof of \crefrange{thm:b}{thm:e}} \label{sec:4}

\begin{proof}
Observe that \(E_0(n) + E_1(n) = E_2(n) + E_3(n) = E(n) = \overline{p_o}(n)\). Thus for \cref{thm:b} and \cref{thm:c} it is sufficient to see that \[E_0(n) - E_1(n) = \begin{cases}
    2, \text{ if \(n\) is a square}, \\
    0, \text{ otherwise};
\end{cases}\]
and for \cref{thm:d} and \cref{thm:e} it is sufficient to see that \[E_2(n) - E_3(n) = \begin{cases}
    2 \cdot (-1)^n, \text{ if \(n\) is a square}, \\
    0, \text{ otherwise}.
\end{cases}\]

Fix any \(\lambda = (\lambda_e, \lambda_g, \lambda_b) \in \mathcal{E}(n)\) and consider its two-modular diagram. Let us agree, that if \(\lambda_e\), \(\lambda_g\) or \(\lambda_b\) is empty, we say that its cardinality is zero.

Let us define a transformation from \(\mathcal{E}(n)\) onto itself:
\begin{enumerate}[label=\arabic*.]
    \item If there are two triangles adjoined, forming a square with diagonal line, then delete the diagonal line.

    \item The next step of the transformation depends on the cardinality of \(\lambda_g\) and \(\lambda_b\).
    
    \begin{enumerate}[label*=\arabic*.]
        \item If \(|\lambda_g| > |\lambda_b|\), then let \(l\) denote the length of the longest column containing only squares. If \(2l\) is strictly greater than the maximum part in \(\lambda_e\) (let us agree that \(\max{\lambda_e} = 0\), if \(\lambda_e = \emptyset\)), then remove that longest column and include \(2l\) as a part in \(\lambda_e\). Otherwise, place the \(\max{\lambda_e}\) part into the diagram as a column to the left of the longest column.

        If \(l = 0\) and \(\lambda_e \neq \emptyset\), then the \(\max{\lambda_e}\) part should be added to the far right of the diagram as a column.

        \item If \(|\lambda_g| \leq |\lambda_b|\), then let \(l\) denote the longest row containing only squares. If \(2l\) is strictly greater than \(\max \lambda_e\), then remove that row and include \(2l\) as a part in \(\lambda_e\). Otherwise, place the \(\max{\lambda_e}\) part as a row into the diagram right above the longest row of squares.

        If \(l = 0\) and \(\lambda_e \neq \emptyset\), then the \(\max{\lambda_e}\) part should be added at the bottom of the diagram as a row.
    \end{enumerate}

    \item After that draw back the diagonal line and read off the new partition \(\lambda' = (\lambda_e', \lambda_g', \lambda_b')\).
\end{enumerate}

There are two partitions for which the transformation is not applicable: \(\lambda = ((2k-1)_g, \dots, 3_g, 1_g)\) and \(\lambda = ((2k-1)_b, \dots, 3_b, 1_b)\). Here the mapping does not work, because such partitions do not have any rows or columns containing only squares and \(\lambda_e = \emptyset\). However, for all other partitions the transformation establishes an involution.

Notice that this transformation does change the parity of \(|\lambda_e|\). Moreover, it is not difficult to see that the value \(|\lambda_g| - |\lambda_b|\) is fixed under the transformation, since \(||\lambda_g| - |\lambda_b||\) is the number of unadjoined triangles along the diagonal line, which is fixed under the transformation; thus it does not change parity of \(|\lambda_g| + |\lambda_b|\).

Thus if we ignore the two aforementioned partitions, we establish a one-to-one correspondence between \(E_0(n) \leftrightarrow E_1(n)\) and \(E_2(n) \leftrightarrow E_3(n)\). 

Note that \(\lambda = ((2k-1)_g, \dots, 3_g, 1_g)\) and \(\lambda = ((2k-1)_b, \dots, 3_b, 1_b)\) are in \(\mathcal{E}(k^2)\), since \(1 + 3 + \dots + (2k-1) = k^2\). Hence the theorem is proved.
\end{proof}

It deserves to be mentioned, that the same idea of an almost-involution happens to be in the proof of Euler's pentagonal theorem \cite{andrews}. Moreover, the generalization of two-modular diagram allows one to prove combinatorically the Jacobi triple product identity \cite{jacobi}. 

Let us demonstrate the transformation established in the proof on some examples.
\begin{itemize}
    \item \(\lambda_e = \{12, 6\}; \lambda_g = \{13, 9, 5, 3\}; \lambda_b = \{5, 1\}\) \(\mapsto\) \(\lambda'_e = \{6\}; \lambda'_g = \{15, 11, 7, 5\}\); \(\lambda'_b = \{7, 3\}\).
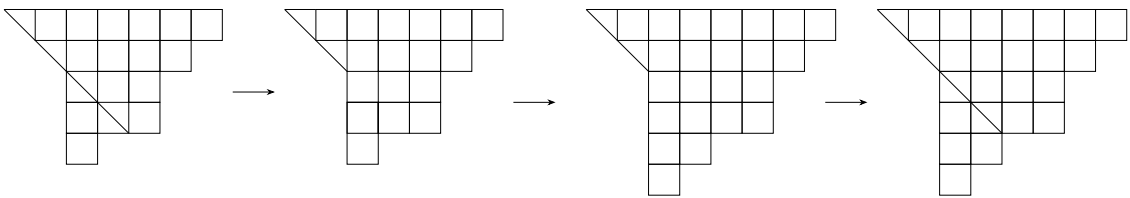
    \begin{figure}[H]
\centering
\resizebox{0.9\textwidth}{!}{%
\begin{circuitikz}
\tikzstyle{every node}=[font=\LARGE]
\draw  (-3,13.25) rectangle (-2.25,12.5);
\draw  (-2.25,13.25) rectangle (-1.5,12.5);
\draw  (-1.5,13.25) rectangle (-0.75,12.5);
\draw  (-0.75,13.25) rectangle (0,12.5);
\draw  (0,13.25) rectangle (0.75,12.5);
\draw  (0.75,13.25) rectangle (1.5,12.5);
\draw  (-2.25,12.5) rectangle (-1.5,11.75);
\draw  (-1.5,12.5) rectangle (-0.75,11.75);
\draw  (-0.75,12.5) rectangle (0,11.75);
\draw  (0,12.5) rectangle (0.75,11.75);
\draw  (-1.5,11.75) rectangle (-0.75,11);
\draw  (-0.75,11.75) rectangle (0,11);
\draw  (-0.75,11) rectangle (0,10.25);
\draw [short] (-3.75,13.25) -- (-0.75,10.25);
\draw [short] (-3.75,13.25) -- (-3,13.25);
\draw [short] (-0.75,10.25) -- (-1.5,10.25);
\draw [short] (-1.5,11) -- (-1.5,10.25);
\draw [short] (-2.25,11.75) -- (-2.25,11);
\draw [short] (-2.25,11) -- (-1.5,11);
\draw  (-2.25,11) rectangle (-1.5,10.25);
\draw  (-2.25,10.25) rectangle (-1.5,9.5);
\draw [->, >=Stealth] (1.75,11.25) -- (2.75,11.25);
\draw  (3.75,13.25) rectangle (4.5,12.5);
\draw  (4.5,13.25) rectangle (5.25,12.5);
\draw  (5.25,13.25) rectangle (6,12.5);
\draw  (6,13.25) rectangle (6.75,12.5);
\draw  (6.75,13.25) rectangle (7.5,12.5);
\draw  (7.5,13.25) rectangle (8.25,12.5);
\draw  (4.5,12.5) rectangle (5.25,11.75);
\draw  (5.25,12.5) rectangle (6,11.75);
\draw  (6,12.5) rectangle (6.75,11.75);
\draw  (6.75,12.5) rectangle (7.5,11.75);
\draw  (5.25,11.75) rectangle (6,11);
\draw  (6,11.75) rectangle (6.75,11);
\draw  (6,11) rectangle (6.75,10.25);
\draw [short] (3,13.25) -- (3.75,13.25);
\draw [short] (5.25,11) -- (5.25,10.25);
\draw [short] (4.5,11.75) -- (4.5,11);
\draw [short] (4.5,11) -- (5.25,11);
\draw  (4.5,11) rectangle (5.25,9.5);
\draw  (4.5,11) rectangle (5.25,10.25);
\draw  (6,10.25) rectangle (5.25,11);
\draw [short] (3,13.25) -- (4.5,11.75);
\draw [->, >=Stealth] (8.5,11) -- (9.5,11);
\draw  (11.75,13.25) rectangle (12.5,12.5);
\draw  (11.75,12.5) rectangle (12.5,11.75);
\draw  (11.75,11.75) rectangle (12.5,11);
\draw  (11.75,11) rectangle (12.5,10.25);
\draw  (11.75,10.25) rectangle (12.5,9.5);
\draw  (11.75,9.5) rectangle (12.5,8.75);
\draw  (12.5,13.25) rectangle (13.25,12.5);
\draw  (12.5,12.5) rectangle (13.25,11.75);
\draw  (12.5,11.75) rectangle (13.25,11);
\draw  (12.5,11) rectangle (13.25,10.25);
\draw  (12.5,10.25) rectangle (13.25,9.5);
\draw  (13.25,13.25) rectangle (14,12.5);
\draw  (13.25,12.5) rectangle (14,11.75);
\draw  (13.25,11.75) rectangle (14,11);
\draw  (13.25,11) rectangle (14,10.25);
\draw  (14,13.25) rectangle (14.75,12.5);
\draw  (14,12.5) rectangle (14.75,11.75);
\draw  (14,11.75) rectangle (14.75,11);
\draw  (14,11) rectangle (14.75,10.25);
\draw  (14.75,13.25) rectangle (15.5,12.5);
\draw  (14.75,12.5) rectangle (15.5,11.75);
\draw  (15.5,13.25) rectangle (16.25,12.5);
\draw  (11.75,13.25) rectangle (11,12.5);
\draw [short] (10.25,13.25) -- (11.75,11.75);
\draw [short] (10.25,13.25) -- (11,13.25);
\draw [->, >=Stealth] (16,11) -- (17,11);
\draw  (18.75,13.25) rectangle (19.5,12.5);
\draw  (18.75,12.5) rectangle (19.5,11.75);
\draw  (18.75,11.75) rectangle (19.5,11);
\draw  (18.75,11) rectangle (19.5,10.25);
\draw  (18.75,10.25) rectangle (19.5,9.5);
\draw  (18.75,9.5) rectangle (19.5,8.75);
\draw  (19.5,13.25) rectangle (20.25,12.5);
\draw  (19.5,12.5) rectangle (20.25,11.75);
\draw  (19.5,11.75) rectangle (20.25,11);
\draw  (19.5,11) rectangle (20.25,10.25);
\draw  (19.5,10.25) rectangle (20.25,9.5);
\draw  (20.25,13.25) rectangle (21,12.5);
\draw  (20.25,12.5) rectangle (21,11.75);
\draw  (20.25,11.75) rectangle (21,11);
\draw  (20.25,11) rectangle (21,10.25);
\draw  (21,13.25) rectangle (21.75,12.5);
\draw  (21,12.5) rectangle (21.75,11.75);
\draw  (21,11.75) rectangle (21.75,11);
\draw  (21,11) rectangle (21.75,10.25);
\draw  (21.75,13.25) rectangle (22.5,12.5);
\draw  (21.75,12.5) rectangle (22.5,11.75);
\draw  (22.5,13.25) rectangle (23.25,12.5);
\draw  (18.75,13.25) rectangle (18,12.5);
\draw [short] (17.25,13.25) -- (18.75,11.75);
\draw [short] (17.25,13.25) -- (18,13.25);
\draw [short] (18.75,11.75) -- (20.25,10.25);
\end{circuitikz}
}%
\label{fig:my_label}
\end{figure}

    \item \(\lambda_e = \emptyset; \lambda_g = \{7, 5, 1\}; \lambda_b = \{5, 3\}\) \(\mapsto\) \(\lambda'_e = \{8\}; \lambda'_g = \{5, 3\}\); \(\lambda'_b = \{5\}\).
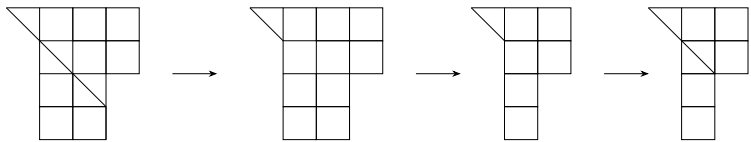
    \begin{figure}[H]
\centering
\resizebox{280pt}{!}{%
\begin{circuitikz}
\tikzstyle{every node}=[font=\LARGE]
\draw  (-2,17) rectangle (-1.25,16.25);
\draw  (-1.25,17) rectangle (-0.5,16.25);
\draw  (-0.5,17) rectangle (0.25,16.25);
\draw  (-1.25,16.25) rectangle (-0.5,15.5);
\draw  (-0.5,16.25) rectangle (0.25,15.5);
\draw [short] (-2.75,17) -- (-2,17);
\draw [short] (-0.5,14.75) -- (-0.5,14);
\draw [short] (-1.25,15.5) -- (-1.25,14.75);
\draw [short] (-1.25,14.75) -- (-0.5,14.75);
\draw  (-1.25,14.75) rectangle (-0.5,14);
\draw [->, >=Stealth] (1,15.5) -- (2,15.5);
\draw [->, >=Stealth] (6.5,15.5) -- (7.5,15.5);

\draw (-2.75,17) to[short] (-0.5,14.75);
\draw (-0.5,15.5) to[short] (-0.5,14.75);
\draw  (-1.25,14.75) rectangle (-2,15.5);
\draw  (-2,14.75) rectangle (-1.25,14);
\draw [short] (-2,16.25) -- (-2,15.5);
\draw  (3.5,17) rectangle (4.25,16.25);
\draw  (4.25,17) rectangle (5,16.25);
\draw  (5,17) rectangle (5.75,16.25);
\draw  (4.25,16.25) rectangle (5,15.5);
\draw  (5,16.25) rectangle (5.75,15.5);
\draw [short] (2.75,17) -- (3.5,17);
\draw [short] (5,14.75) -- (5,14);
\draw [short] (4.25,15.5) -- (4.25,14.75);
\draw [short] (4.25,14.75) -- (5,14.75);
\draw  (4.25,14.75) rectangle (5,14);
\draw (5,15.5) to[short] (5,14.75);
\draw  (4.25,14.75) rectangle (3.5,15.5);
\draw  (3.5,14.75) rectangle (4.25,14);
\draw [short] (3.5,16.25) -- (3.5,15.5);
\draw [short] (2.75,17) -- (3.5,16.25);
\draw  (8.5,17) rectangle (9.25,16.25);
\draw  (8.5,16.25) rectangle (9.25,15.5);
\draw  (8.5,15.5) rectangle (9.25,14.75);
\draw  (8.5,14.75) rectangle (9.25,14);
\draw  (9.25,17) rectangle (10,16.25);
\draw  (9.25,16.25) rectangle (10,15.5);
\draw [short] (7.75,17) -- (8.5,16.25);
\draw [short] (7.75,17) -- (8.5,17);
\draw  (12.5,17) rectangle (13.25,16.25);
\draw  (12.5,16.25) rectangle (13.25,15.5);
\draw  (12.5,15.5) rectangle (13.25,14.75);
\draw  (12.5,14.75) rectangle (13.25,14);
\draw  (13.25,17) rectangle (14,16.25);
\draw  (13.25,16.25) rectangle (14,15.5);
\draw [short] (11.75,17) -- (12.5,16.25);
\draw [short] (11.75,17) -- (12.5,17);
\draw [short] (12.5,16.25) -- (13.25,15.5);
\draw [->, >=Stealth] (10.75,15.5) -- (11.75,15.5);
\end{circuitikz}
}%

\label{fig:my_label}
\end{figure}

    \newpage
    
    \item \(\lambda_e = \{12, 2\}; \lambda_g = \{7, 3\}; \lambda_b = \{11, 7, 5, 3\}\) \(\mapsto\) \(\lambda'_e = \{2\}; \lambda'_g = \{7, 5, 1\}\); \(\lambda'_b = \{13, 9, 7, 5, 1\}\).
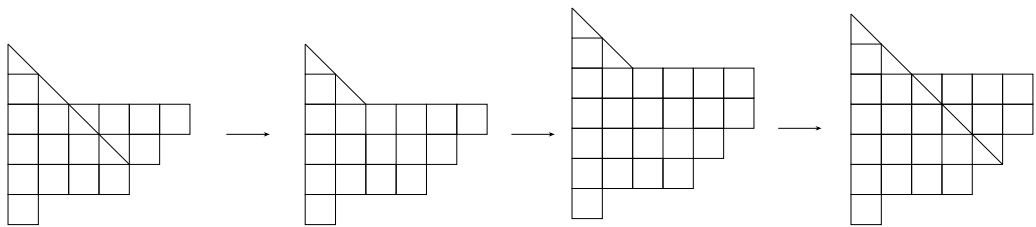
    \begin{figure}[H]
\centering
\resizebox{385pt}{!}{%
\begin{circuitikz}
\tikzstyle{every node}=[font=\LARGE]
\draw  (-15.75,27.75) rectangle (-14.5,26.5);
\draw  (-15.75,26.5) rectangle (-14.5,25.25);
\draw  (-15.75,25.25) rectangle (-14.5,24);
\draw  (-15.75,24) rectangle (-14.5,22.75);
\draw  (-15.75,22.75) rectangle (-14.5,21.5);
\draw  (-14.5,26.5) rectangle (-13.25,25.25);
\draw  (-14.5,25.25) rectangle (-13.25,24);
\draw  (-14.5,24) rectangle (-13.25,22.75);
\draw  (-13.25,25.25) rectangle (-12,24);
\draw  (-13.25,24) rectangle (-12,22.75);
\draw  (-12,24) rectangle (-10.75,22.75);
\draw [short] (-10.75,24) -- (-15.75,29);
\draw [short] (-15.75,27.75) -- (-15.75,29);
\draw  (-10.75,24) rectangle (-9.5,25.25);
\draw  (-9.5,25.25) rectangle (-10.75,26.5);
\draw  (-9.5,26.5) rectangle (-8.25,25.25);
\draw  (-10.75,26.5) rectangle (-12,25.25);
\draw [short] (-13.25,26.5) -- (-12,26.5);
\draw [->, >=Stealth] (-6.75,25.25) -- (-5,25.25);
\draw  (-3.5,27.75) rectangle (-2.25,26.5);
\draw  (-3.5,26.5) rectangle (-2.25,25.25);
\draw  (-3.5,25.25) rectangle (-2.25,24);
\draw  (-3.5,24) rectangle (-2.25,22.75);
\draw  (-3.5,22.75) rectangle (-2.25,21.5);
\draw  (-2.25,26.5) rectangle (-1,25.25);
\draw  (-2.25,25.25) rectangle (-1,24);
\draw  (-2.25,24) rectangle (-1,22.75);
\draw  (-1,25.25) rectangle (0.25,24);
\draw  (-1,24) rectangle (0.25,22.75);
\draw  (0.25,24) rectangle (1.5,22.75);
\draw [short] (-3.5,27.75) -- (-3.5,29);
\draw  (1.5,24) rectangle (2.75,25.25);
\draw  (2.75,25.25) rectangle (1.5,26.5);
\draw  (2.75,26.5) rectangle (4,25.25);
\draw  (1.5,26.5) rectangle (0.25,25.25);
\draw [short] (-1,26.5) -- (0.25,26.5);
\draw [short] (-3.5,29) -- (-1,26.5);
\draw [->, >=Stealth] (5.0,25.25) -- (6.75,25.25);
\draw  (7.5,26.75) rectangle (8.75,25.5);
\draw  (7.5,25.5) rectangle (8.75,24.25);
\draw  (7.5,24.25) rectangle (8.75,23);
\draw  (7.5,23) rectangle (8.75,21.75);
\draw  (8.75,26.75) rectangle (10,25.5);
\draw  (8.75,25.5) rectangle (10,24.25);
\draw  (8.75,24.25) rectangle (10,23);
\draw  (10,25.5) rectangle (11.25,24.25);
\draw  (10,24.25) rectangle (11.25,23);
\draw  (11.25,24.25) rectangle (12.5,23);
\draw  (12.5,24.25) rectangle (13.75,25.5);
\draw  (13.75,25.5) rectangle (12.5,26.75);
\draw  (13.75,26.75) rectangle (15,25.5);
\draw  (12.5,26.75) rectangle (11.25,25.5);
\draw [->, >=Stealth] (16.0,25.5) -- (17.75,25.5);
\draw  (10,26.75) rectangle (11.25,25.5);
\draw  (7.5,26.75) rectangle (8.75,28);
\draw  (8.75,28) rectangle (10,26.75);
\draw  (10,28) rectangle (11.25,26.75);
\draw  (11.25,28) rectangle (12.5,26.75);
\draw  (12.5,28) rectangle (13.75,26.75);
\draw  (13.75,28) rectangle (15,26.75);
\draw  (7.5,28) rectangle (8.75,29.25);
\draw [short] (10,28) -- (7.5,30.5);
\draw [short] (7.5,29.25) -- (7.5,30.5);
\draw  (19,26.5) rectangle (20.25,25.25);
\draw  (19,25.25) rectangle (20.25,24);
\draw  (19,24) rectangle (20.25,22.75);
\draw  (19,22.75) rectangle (20.25,21.5);
\draw  (20.25,26.5) rectangle (21.5,25.25);
\draw  (20.25,25.25) rectangle (21.5,24);
\draw  (20.25,24) rectangle (21.5,22.75);
\draw  (21.5,25.25) rectangle (22.75,24);
\draw  (21.5,24) rectangle (22.75,22.75);
\draw  (22.75,24) rectangle (24,22.75);
\draw  (24,24) rectangle (25.25,25.25);
\draw  (25.25,25.25) rectangle (24,26.5);
\draw  (25.25,26.5) rectangle (26.5,25.25);
\draw  (24,26.5) rectangle (22.75,25.25);
\draw  (21.5,26.5) rectangle (22.75,25.25);
\draw  (19,26.5) rectangle (20.25,27.75);
\draw  (20.25,27.75) rectangle (21.5,26.5);
\draw  (21.5,27.75) rectangle (22.75,26.5);
\draw  (22.75,27.75) rectangle (24,26.5);
\draw  (24,27.75) rectangle (25.25,26.5);
\draw  (25.25,27.75) rectangle (26.5,26.5);
\draw  (19,27.75) rectangle (20.25,29);
\draw [short] (21.5,27.75) -- (19,30.25);
\draw [short] (19,29) -- (19,30.25);
\draw [short] (21.5,27.75) -- (25.25,24);
\end{circuitikz}
}%

\label{fig:my_label}
\end{figure}
    
    \item \(\lambda_e = \{4\}; \lambda_g = \{11, 7, 3, 1\}; \lambda_b = \{9, 7, 3, 1\}\) \(\mapsto\) \(\lambda'_e = \{12, 4\}; \lambda'_g = \{9, 5, 3\}\); \(\lambda'_b = \{7, 5, 1\}\).
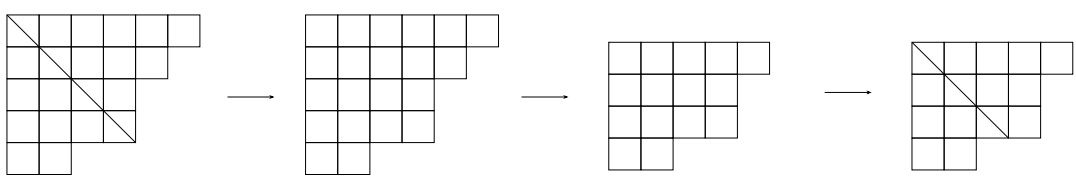
    \begin{figure}[H]
\centering
\resizebox{400pt}{!}{%
\begin{circuitikz}
\tikzstyle{every node}=[font=\LARGE]
\draw  (-16.5,33.5) rectangle (-14.75,31.75);
\draw  (-14.75,33.5) rectangle (-13,31.75);
\draw  (-13,33.5) rectangle (-11.25,31.75);
\draw  (-11.25,33.5) rectangle (-9.5,31.75);
\draw  (-9.5,33.5) rectangle (-7.75,31.75);
\draw  (-14.75,31.75) rectangle (-13,30);
\draw  (-13,31.75) rectangle (-11.25,30);
\draw  (-11.25,31.75) rectangle (-9.5,30);
\draw  (-13,30) rectangle (-11.25,28.25);
\draw  (-13,28.25) rectangle (-11.25,26.5);
\draw  (-13,28.25) rectangle (-14.75,26.5);
\draw  (-14.75,28.25) rectangle (-16.5,26.5);
\draw  (-16.5,28.25) rectangle (-18.25,26.5);
\draw  (-18.25,26.5) rectangle (-16.5,24.75);
\draw  (-16.5,26.5) rectangle (-14.75,24.75);
\draw  (-14.75,30) rectangle (-13,28.25);
\draw  (-14.75,28.25) rectangle (-16.5,30);
\draw  (-16.5,30) rectangle (-18.25,28.25);
\draw  (-14.75,31.75) rectangle (-16.5,30);
\draw  (-16.5,31.75) rectangle (-18.25,30);
\draw  (-16.5,33.5) rectangle (-18.25,31.75);
\draw [short] (-18.25,33.5) -- (-11.25,26.5);
\draw [->, >=Stealth] (-6.25,29) -- (-3.75,29);
\draw  (-0.25,33.5) rectangle (1.5,31.75);
\draw  (1.5,33.5) rectangle (3.25,31.75);
\draw  (3.25,33.5) rectangle (5,31.75);
\draw  (5,33.5) rectangle (6.75,31.75);
\draw  (6.75,33.5) rectangle (8.5,31.75);
\draw  (1.5,31.75) rectangle (3.25,30);
\draw  (3.25,31.75) rectangle (5,30);
\draw  (5,31.75) rectangle (6.75,30);
\draw  (3.25,30) rectangle (5,28.25);
\draw  (3.25,28.25) rectangle (5,26.5);
\draw  (3.25,28.25) rectangle (1.5,26.5);
\draw  (1.5,28.25) rectangle (-0.25,26.5);
\draw  (-0.25,28.25) rectangle (-2,26.5);
\draw  (-2,26.5) rectangle (-0.25,24.75);
\draw  (-0.25,26.5) rectangle (1.5,24.75);
\draw  (1.5,30) rectangle (3.25,28.25);
\draw  (1.5,28.25) rectangle (-0.25,30);
\draw  (-0.25,30) rectangle (-2,28.25);
\draw  (1.5,31.75) rectangle (-0.25,30);
\draw  (-0.25,31.75) rectangle (-2,30);
\draw  (-0.25,33.5) rectangle (-2,31.75);
\draw [->, >=Stealth] (9.75,29) -- (12.25,29);
\draw  (18,32) rectangle (19.75,30.25);
\draw  (19.75,32) rectangle (21.5,30.25);
\draw  (21.5,32) rectangle (23.25,30.25);
\draw  (19.75,30.25) rectangle (21.5,28.5);
\draw  (19.75,28.5) rectangle (21.5,26.75);
\draw  (19.75,28.5) rectangle (18,26.75);
\draw  (18,28.5) rectangle (16.25,26.75);
\draw  (16.25,28.5) rectangle (14.5,26.75);
\draw  (14.5,26.75) rectangle (16.25,25);
\draw  (16.25,26.75) rectangle (18,25);
\draw  (18,30.25) rectangle (19.75,28.5);
\draw  (18,28.5) rectangle (16.25,30.25);
\draw  (16.25,30.25) rectangle (14.5,28.5);
\draw  (18,32) rectangle (16.25,30.25);
\draw  (16.25,32) rectangle (14.5,30.25);
\draw [->, >=Stealth] (26.25,29.25) -- (28.75,29.25);
\draw  (34.5,32) rectangle (36.25,30.25);
\draw  (36.25,32) rectangle (38,30.25);
\draw  (38,32) rectangle (39.75,30.25);
\draw  (36.25,30.25) rectangle (38,28.5);
\draw  (36.25,28.5) rectangle (38,26.75);
\draw  (36.25,28.5) rectangle (34.5,26.75);
\draw  (34.5,28.5) rectangle (32.75,26.75);
\draw  (32.75,28.5) rectangle (31,26.75);
\draw  (31,26.75) rectangle (32.75,25);
\draw  (32.75,26.75) rectangle (34.5,25);
\draw  (34.5,30.25) rectangle (36.25,28.5);
\draw  (34.5,28.5) rectangle (32.75,30.25);
\draw  (32.75,30.25) rectangle (31,28.5);
\draw  (34.5,32) rectangle (32.75,30.25);
\draw  (32.75,32) rectangle (31,30.25);
\draw [short] (31,32) -- (36.25,26.75);
\end{circuitikz}
}%

\label{fig:my_label}
\end{figure}
    
    \item \(\lambda_e = \{2\}; \lambda_g = \{3, 1\}; \lambda_b = \emptyset\) \(\mapsto\) \(\lambda'_e = \emptyset; \lambda'_g = \{5, 1\}\); \(\lambda'_b = \emptyset\).
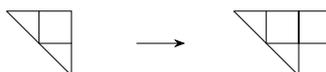
    \begin{figure}[H]
\centering
\resizebox{125pt}{!}{%
\begin{circuitikz}
\tikzstyle{every node}=[font=\LARGE]
\draw  (4.75,10) rectangle (5.25,9.5);
\draw [short] (4.75,10) -- (4.25,10);
\draw [short] (4.25,10) -- (5.25,9);
\draw [short] (5.25,9) -- (5.25,9.5);
\draw  (8.25,10) rectangle (8.75,9.5);
\draw  (8.75,10) rectangle (9.25,9.5);
\draw [short] (8.25,10) -- (7.75,10);
\draw [short] (7.75,10) -- (8.75,9);
\draw [short] (8.75, 9) -- (8.75, 10);
\draw [->, >=Stealth] (6.25,9.5) -- (7,9.5);
\end{circuitikz}
}%

\label{fig:my_label}
\end{figure}

\end{itemize}

\section{Acknowledgments}

I would like to thank my advisor Eric Mortenson for suggesting the problems in the paper of Andrews and El Bachraoui.

This work was performed at the Saint Petersburg Leonhard Euler International Mathematical Institute and supported by the Ministry of Science and Higher Education of the Russian Federation (agreement no. 075-15-2025-343).

\section{Data availability statement}

The authors declare that the data supporting the findings of this study are available within the paper.

\bibliographystyle{plaindin}
\bibliography{sample}

@misc{andrews2024twocolorpartitionsoddsmallest,
      title={On two-color partitions with odd smallest part}, 
      author={George E. Andrews and Mohamed {El Bachraoui}},
      year={2024},
      eprint={2410.14190},
      archivePrefix={arXiv},
      primaryClass={math.NT},
      url={https://arxiv.org/abs/2410.14190}, 
}

@article{articlecombinatoricalproofs,
    author = {Chen, Dandan and Zou, Ziyin},
    year = {2025},
    month = {04},
    pages = {1-5},
    title = {\emph{Combinatorial proofs for two-colour partitions}},
    journal = {Bulletin of the Australian Mathematical Society},
    doi = {10.1017/S0004972725000206}
    }

@article{jacobi,
    author = {Kolitsch, Louis W and Kolitsch, Stephanie},
    title = {\emph{A combinatorial proof of {J}acobi’s triple product identity}},
    journal = {The Ramanujan Journal},
    year = {2018},
    pages = {483 - 489},
    doi = {10.1007/s11139-016-9854-5}
}

@book{andrews,
author = {George E. Andrews},
year = {1976},
title = {The Theory of Partitions},
publisher = {Cambridge University Press}
}

\end{document}